\documentclass{birkmult}
\usepackage{graphicx,amssymb}
\usepackage{amsmath,amssymb,amscd}
\input xy
\xyoption{all}
\usepackage[all]{xy}
\usepackage{hyperref}

\newcommand{\ra}{\rightarrow}

\newcommand{\cO}{\mathcal{O}}
\newcommand{\cT}{\mathcal{T}}

\theoremstyle{plain}
\newtheorem{theorem}{Theorem}[section]
\newtheorem{lem}[theorem]{Lemma}
\newtheorem{prop}[theorem]{Proposition}
\newtheorem{cor}[theorem]{Corollary}
\newtheorem{defin}[theorem]{Definition}
\newtheorem{rem}[theorem]{Remark}

\numberwithin{equation}{section}

\begin{document}
\title[Clifford indices]{Generation of vector bundles computing Clifford indices}

\author{H. Lange}
\author{P. E. Newstead}

\address{H. Lange\\Department Mathematik\\
              Universit\"at Erlangen-N\"urnberg\\
              Bismarckstra\ss e $1\frac{ 1}{2}$\\
              D-$91054$ Erlangen\\
              Germany}
              \email{lange@mi.uni-erlangen.de}
\address{P.E. Newstead\\Department of Mathematical Sciences\\
              University of Liverpool\\
              Peach Street, Liverpool L69 7ZL, UK}
\email{newstead@liv.ac.uk}

\thanks{Both authors are members of the research group VBAC (Vector Bundles on Algebraic Curves). The second author 
would like to thank the Department Mathematik der Universit\"at 
         Erlangen-N\"urnberg for its hospitality}

\begin{abstract}
Clifford indices for semistable vector bundles on a smooth projective curve of genus at least 4 were defined
in a previous paper of the authors. The present paper studies bundles which compute these Clifford indices.
We show that under certain conditions on the curve all such bundles and their Serre duals are generated. 
\end{abstract}
\maketitle

\section{Introduction}

Let $C$ be a smooth projective curve of genus $g \geq 4$ defined over an algebraically closed field of characteristic zero. 
The classical Clifford index of $C$ has played a major role in describing the geometry of $C$. Some years ago E. Ballico 
\cite{b} proposed a number of possible generalisations of the Clifford index to vector bundles but did not develop this 
idea to any extent. A common feature of most of these definitions was that the bundles considered were required to be 
spanned (or, as we shall say, generated), that is, that the evaluation map $H^0(E)\otimes {\mathcal O}_C\to E$ is 
surjective. This is a very natural condition as generated bundles are precisely those we need in order to define 
morphisms from $C$ to a Grassmannian. In fact, one of Ballico's definitions restricted further to primitive 
vector bundles (in other words, bundles $E$ for which $E$ and $E^*\otimes K$ are both generated); this generalises 
a well known concept for line bundles (see \cite{ckm}). Ballico used stable bundles in three of his definitions, 
but never semistable bundles. 

In a previous paper \cite{cl}, we proposed two definitions of Clifford index for semistable bundles in the 
following way. First we define, for any vector bundle $E$ of rank $n$ and degree $d$, 
$$
\gamma(E) := \frac{1}{n} \left(d - 2(h^0(E) -n)\right) = \mu(E) -2\frac{h^0(E)}{n} + 2,
$$
where $\mu(E)=\frac{d}n$.
Ballico defines $\text{Cliff}(E)$ in the same way but without the scaling factor $\frac{1}{n}$, 
which we include since it makes it easier to compare Clifford indices for different ranks. Then 
the Clifford indices $\gamma_n$ and $\gamma_n'$ are defined by 
$$
\gamma_n := \min_{E} \left\{ \gamma(E) \left|  
\begin{array}{c}   E \;\mbox{semistable of rank}\; n \\
h^0(E) \geq n+1,\; \mu(E) \leq g-1
\end{array} \right\} \right.
$$
and
$$
\gamma_n' := \min_{E} \left\{ \gamma(E) \;\left| 
\begin{array}{c} E \;\mbox{semistable of rank}\; n \\
h^0(E) \geq 2n,\; \mu(E) \leq g-1
\end{array} \right\}. \right.
$$
In particular we did not make any assumption of generation or primitivity. It may be noted 
that $\gamma_1=\gamma_1'$ is the usual Clifford index. We say that a bundle $E$ {\em contributes} 
to $\gamma_n$ if it is semistable of rank $n$ with $\mu(E) \leq g-1$ and $h^0(E)\ge n+1$ and that $E$ 
{\em computes} $\gamma_n$ if in addition $\gamma(E)=\gamma_n$. Similar definitions are made for $\gamma'_n$. 

Our first object in this paper is to prove that, under certain conditions, every bundle computing 
$\gamma_n'$ or $\gamma_n$ is primitive, so that we do not need to assume primitivity or generation 
in the definition. More precisely, we show\\

\noindent{\bf Theorem 2.4.}
\begin{em} Suppose $\gamma'_r\ge\gamma_n'$ for all $r\le n$. If $E$ is a semistable bundle computing $\gamma'_n$, 
then $E$ is primitive.
\end{em}\\

\noindent{\bf Theorem 2.5.}
\begin{em} Suppose $\gamma_r\ge\gamma_n$ for all $r\le n$. If $E$ is a semistable bundle computing $\gamma_n$, 
then $E$ is primitive.
\end{em}\\

The hypotheses of these theorems are satisfied for $n=2$ (Corollary 2.6) and that of Theorem 2.5 is satisfied 
whenever $n\ge g-3$ (Corollary 2.7). For $n=3$, the hypothesis of Theorem 2.4 holds when $\gamma_1\le4$ and 
when $C$ is a smooth plane curve of degree $\ge5$. In fact, we know of no examples where the hypothesis of Theorem 2.4 fails. For $n=3$, the 
hypothesis of Theorem 2.5 holds for a general curve (Corollary 2.10), but we shall indicate an example 
(for $C$ a smooth plane curve) for which neither the hypothesis nor the conclusion holds.

If the hypotheses of the theorems fail to hold, one may still ask whether a bundle $E$ computing $\gamma_n$ or 
$\gamma_n'$ is generically generated (that is, the evaluation map of $E$ has torsion cokernel). In section 3 
we will prove that this is so in rank $3$:\\

\noindent{\bf Theorem 3.3.}
\begin{em} Any semistable bundle of rank $3$ computing either $\gamma'_3$ or $\gamma_3$ is generically generated.
\end{em}\\

\noindent{\bf Notation.} Throughout the paper $C$ will denote a smooth projective curve of genus $g \geq 4$ defined 
over an algebraically closed field of characteristic zero. For a vector bundle $G$ on $C$, the rank and 
degree of $G$ will be denoted by $r_G$ and $d_G$ respectively. The {\em gonality sequence} 
$d_1,d_2,\ldots,d_r,\ldots$ of $C$ is defined by 
$$
d_r := \min \{ d_L \;|\; L \; \mbox{a line bundle on} \; C \; \mbox{with} \; h^0(L) \geq r +1\}.
$$
We have always $d_{r+s}\le d_r+d_s$ and in particular $d_n\le nd_1$ for all $n$ (see \cite[Section 4]{cl}).

\section{Generation of bundles computing $\gamma'_n$}

Fix positive integers $d,n$ and $s$. 

\begin{lem} \label{lem2.1}
Suppose there exists a bundle $F$ of rank $r_F \leq n$  which satisfies either $h^0(F) > (\frac{s}{n} +1)r_F$ 
and every subbundle of $F$ has slope $\leq \frac{d}{n}$
or $h^0(F) \geq (\frac{s}{n} +1)r_F$ and every subbundle of $F$ has slope $< \frac{d}{n}$. Then
$$
\frac{d-2s}{n} > \min_{G \in \cT} \gamma(G)
$$
where 
$$
\cT = \left\{ G \; \mbox{semistable} \;| \; r_G \le n,\; \frac{d_G}{r_G} \leq \frac{d}{n}, \; 
h^0(G) \geq \left( \frac{s}{n} + 1 \right) r_G \right\}.
$$
\end{lem}

\begin{proof}
If $F$ is semistable, then either
$$
\gamma(F) = \frac{1}{r_F}\left( d_F - 2(h^0(F) -r_F) \right)  < \frac{d_F}{r_F} - \frac{2}{r_F} \frac{s}{n}r_F \leq \frac{d}{n} - \frac{2s}{n}.
$$
or 
$$
\gamma(F) \leq \frac{d_F}{r_F} - \frac{2}{r_F} \frac{s}{n}r_F < \frac{d}{n} - \frac{2s}{n}.
$$
Moreover, $F \in \cT$ and the result follows.

If $F$ is not semistable, we will apply induction on the rank. So assume that the assertion is true for any 
bundle of smaller rank satisfying the assumptions.

We can write
\begin{equation} \label{eq2.1}
0 \ra M \ra F \ra N \ra 0
\end{equation}
with $M$ the first term in the Harder-Narasimhan filtration. 
                                                                                                                                                      
By definition of the Harder-Narasimham filtration every subbundle of $N$ has slope $< \frac{d_M}{r_M} \leq \frac{d}{n}$.
If $h^0(N) \geq  (\frac{s}{n} +1)r_N$, the induction hypothesis can be applied to $N$ to give 
$$
\frac{d-2s}{n} > \min_{G \in \cT} \gamma(G).
$$

If $h^0(N) < (\frac{s}{n} + 1)r_N$, then 
$$
h^0(M) > \left( \frac{s}{n} + 1 \right) r_F - \left( \frac{s}{n}+1 \right) r_N = \left( \frac{s}{n} +1 \right) r_M.
$$ 
But $M$ is semistable, so we get the conclusion as in the first part of the proof. 
\end{proof}

\begin{prop} \label{prop2.2} 
If $E$ is a semistable bundle of rank $n$ and degree $d$ with $h^0(E) = n+s,\; s \geq 1$ and such that
$$
\gamma(E) \leq \min_{G \in \cT} \gamma(G),
$$
then $E$ is generated.
\end{prop} 

\begin{proof}
Suppose first that $F$ is a proper subbundle of $E$ of rank $r_F$ and degree $d_F$ with $h^0(F) = h^0(E)$. Then 
$$
h^0(F) = n + s > r_F + \frac{sr_F}{n} = \left( \frac{s}{n} + 1 \right) r_F.
$$
Since also $\frac{d_F}{r_F} \leq \frac{d}{n}$ by semistability of $E$, Lemma \ref{lem2.1} applies to give 
$$
\gamma(E) = \frac{d -2s}{n} > \min_{G \in \cT} \gamma(G),
$$
contradicting the hypotheses.

Now suppose $F$ is a subsheaf of $E$ of rank $r_F = n$ and degree $d_F < d$ with $h^0(F) = h^0(E)$. If $F$ is 
semistable, then $F \in \cT$ and $\gamma(F) < \gamma(E)$, a contradiction. 

If $F$ is not 
semistable, then there is an exact sequence \eqref{eq2.1}.
If $h^0(M) > (\frac{s}{n} + 1)r_M$, then $M \in \cT$ and $\gamma(M) < \gamma(E)$ gives the same contradiction.

If $h^0(M) \leq (\frac{s}{n} + 1)r_M$, then $h^0(N) \geq (\frac{s}{n} + 1) r_N$ as above. 
By definition of the Harder-Narasimhan filtration, every subbundle of $N$ has slope $< \frac{d}{n}$. 
Lemma \ref{lem2.1} applied to $N$ now gives a contradiction.
\end{proof}

Recall that a line bundle $L$ on $C$ is called primitive (see \cite[(1.1)]{ckm}), if $L$ and $L ^* \otimes K$ are both generated.
In \cite[Definition 1.2 and Remark 1.3]{b} Ballico generalised this notion to vector bundles. 

\begin{defin}
{\em A vector bundle $E$ on $C$ is called {\em primitive} if $E$ and $E^* \otimes K$ are both generated.}
\end{defin}

Note that $E$ is primitive if and only if $E$ is generated and  for any vector bundle $F$ containing $E$ with $F/E$ 
a torsion sheaf of length 1 we have $h^0(F) = h^0(E)$.

\begin{theorem} \label{thm2.3}
Suppose $\gamma'_r \geq \gamma'_n$ for all $r \leq n$. If $E$ is a semistable bundle
computing $\gamma'_n$, then $E$ is primitive.
\end{theorem}

\begin{proof}
Let $d_E = d$ and $h^0(E) = n + s$ where now $s \geq n$. Also $\frac{d}{n} \leq g-1$. It follows that, if $G \in \cT$, then $G$ contributes
to $\gamma'_{r_G}$. So
$$
\min_{G \in \cT} \gamma(G) \geq \min_{r \leq n} \gamma'_r = \gamma'_n = \gamma(E).
$$
Now Proposition \ref{prop2.2} implies that $E$ is generated.

Now suppose $d = d_{E^* \otimes K}$ and $ h^0(E^* \otimes K) = n+s$. Then $d\ge n(g-1)$ and 
$$
h^0(E) = h^1(E^* \otimes K) = n+s-d+n(g-1).
$$
The condition for $E$ to contribute to $\gamma'_n$ is that 
\begin{equation} \label{eq2.2}
s - d+n(g-1) \geq n.
\end{equation}
So, if $G \in \cT$ and $\frac{d_G}{r_G} \leq g-1$, then $h^0(G) \geq (\frac{s}{n} + 1)r_G \geq 2r_G$. 
Hence $G$ contributes to $\gamma'_{r_G}$. 

If $\frac{d_G}{r_G} > g-1$, then using \eqref{eq2.2} we have
\begin{eqnarray*}
h^0(G^* \otimes K) & = & h^0(G) - d_G + r_G(g-1) \\
& \geq & \left(\frac{s}{n} + 1 \right)r_G -d_G + r_G(g-1) \\
& \geq & \left(\frac{d}{n} - \frac{d_G}{r_G} + 1 \right)r_G -(g-2)r_G + r_G(g-1) \geq 2r_G.
\end{eqnarray*}
So $G^* \otimes K$ contributes to $\gamma'_{r_G}$. Hence, for any $G \in \cT$, 
$$
\gamma(G) = \gamma(G^* \otimes K) \geq \gamma'_{r_G} \geq \gamma'_n = \gamma(E) = \gamma(E^* \otimes K).
$$
Now Proposition \ref{prop2.2} applied to $E^* \otimes K$ shows that $E^* \otimes K$ is generated.
\end{proof}

\begin{theorem} \label{thm2.4}
Suppose $\gamma_r \geq \gamma_n$ for all $r \leq n$. If $E$ is a semistable bundle
computing $\gamma_n$, then $E$ is primitive.
\end{theorem}

\begin{proof}
Let $h^0(E) = n+s$, where $s \geq 1$. The proof proceeds exactly as the proof of the previous theorem using 
$s -d +n(g-1) \geq 1$ instead of \eqref{eq2.2}. 
\end{proof}

\begin{cor} \label{cor2.5}
If $E$ is semistable of rank $2$ and computes either $\gamma_2$ or $\gamma'_2$, then $E$ is primitive.
\end{cor}
\begin{proof}
This follows from the fact that $\gamma_2 \leq \gamma'_2 \leq \gamma_1$ (see \cite[Lemma 2.2]{cl}).
\end{proof}

\begin{cor} \label{cor2.6}
If $E$ is semistable of rank $n \geq g-3$ and computes $\gamma_n$, then $E$ is primitive.
\end{cor}

\begin{proof}
From \cite[Theorems 3.6 and 4.21]{cl} we know that 
$$
\gamma_r \geq \gamma_n \quad \mbox{for all} \quad r \leq n
$$ 
if $n \geq g-3$. So the corollary is a 
consequence of Theorem \ref{thm2.4}.
\end{proof}

\begin{rem} {\em (1) For $ r_E = 3$ the hypothesis of Theorem 2.4 reduces to $\gamma'_3 \leq \gamma'_2$ according to \cite[Lemma 2.2]{cl}. 
This is valid whenever $\gamma'_2 = \gamma_1$ (again by \cite[Lemma 2.2]{cl}) and this is true if $\gamma_1 \leq 4$ 
by \cite[Proposition 3.8]{cl} and also for any smooth plane curve of degree $\geq 5$ by \cite[Proposition 8.1]{cl}.

(2) For $r_E = 4$ the hypothesis reduces to $\gamma'_4 \leq \gamma'_3$ according to \cite[Lemma 2.2]{cl}. 

(3) We know of no examples for which the hypothesis fails.
}
\end{rem}

\begin{prop}
Suppose $\frac{1}{3}(d_3-2) \leq \frac{1}{2}(d_2-2)$. If $E$ is semistable of rank $3$ and computes $\gamma_3$,
then $E$ is primitive.
\end{prop}

\begin{proof}
By Theorem \ref{thm2.4} and \cite[Lemma 2.2]{cl} it is sufficient to show that $\gamma_3 \leq \gamma_2$. 
Note that under the stated hypothesis, $\frac{d_3}{3} < \frac{d_2}{2}$. So
according to \cite[Theorem 6.1 and Corollary 5.3]{cl},
$$
\gamma_3 = \min \left\{ \gamma'_3, \frac{1}{3}(d_3-2) \right\} \quad \mbox{and} \quad 
\gamma_2 = \min \left\{ \gamma_1, \frac{1}{2}(d_2 -2) \right\}.
$$ 
Since $\gamma'_3 \leq \gamma_1$ by \cite[Lemma 2.2]{cl}, it follows from the hypothesis that $\gamma_3 \leq \gamma_2$.
\end{proof}

\begin{cor}
Suppose $C$ is general and $E$ is semistable of rank $3$ computing $\gamma_3$. Then $E$ is primitive.
\end{cor}

\begin{proof}
For a general curve the values of $d_r$ are given by 
\begin{equation} \label{dr}
d_r = r + g - \left[ \frac{g}{r+1} \right]
\end{equation}
(see \cite[Remark 4.4(c)]{cl}).
The condition $\frac{1}{3}(d_3 - 2) \leq \frac{1}{2}(d_2 - 2)$ can be restated as
$$
2\left( 1 + g - \left[ \frac{g}{4} \right] \right) \leq 3 \left( g - \left[ \frac{g}{3} \right] \right),
$$
which is easy to verify for any $g \geq 4$. 
\end{proof}

\begin{rem}
{\em The condition $\frac{1}{3}(d_3 -2) \leq \frac{1}{2}(d_2 - 2)$ is slightly stronger than the hypothesis 
$\frac{d_3}{3} \leq \frac{d_2}{2}$ used repeatedly in \cite[Section 4]{cl}.
}
\end{rem}

\section{Generic generation for rank-3 bundles}

We begin with two lemmas which together formalise an argument in \cite[Section 8]{cl}. 
The first of these is certainly well known.

\begin{lem} \label{lem3.1}
Let $F$ be a stable bundle of rank $2$ and odd degree $d_F$ and let $L$ be a line bundle of degree $\frac{d_F + 1}{2}$.
For any non-trivial extension
\begin{equation} \label{eq3.1}
0 \ra F \ra E' \ra L \ra 0,
\end{equation}
$E'$ is stable.
\end{lem}

\begin{proof}
Suppose first that $M$ is a subbundle of $E'$ of rank 1. 
If $M \subset F$, then by stability of $F$, $d_M < \mu(F) < \mu(E')$. If the map $M \ra L$ is non-zero, then,
since the extension is non-trivial,
$$
d_M \leq \frac{d_F-1}{2} < \mu(E').
$$

Suppose now that $G \subset E'$ is a subbundle of rank 2. We may assume that $G \neq F$. So the map $G \ra L$ is non-zero. 
Then consider the diagram with exact rows,
$$
\xymatrix{
0 \ar[r] & F' \ar[r] \ar@{^{(}->}[d] & G \ar[r] \ar@{^{(}->}[d]  & L' \ar[r] \ar@{^{(}->}[d] & 0\\
0 \ar[r] & F \ar[r] & E' \ar[r] & L \ar[r] & 0.
}
$$
So $d_{F'} < \frac{d_F}{2}$ and $ d_{L'} \leq \frac{d_F+1}{2}$ which implies $d_G < d_F + \frac{1}{2}$ and hence
$\mu(G) \leq \frac{d_F}{2} < \mu(E')$. 
\end{proof}

\begin{lem} \label{lem3.2}
Suppose $F$ is a stable bundle of rank $2$ and odd degree $d_F$ with $1 < d_F \leq 4g-4$. 
For any line bundle $L$ of degree $\frac{d_F + 1}{2}$ with $h^0(L) \geq 1$, there is a non-trivial extension
\eqref{eq3.1} such that a non-zero section of $L$ lifts to $E'$.
\end{lem}

\begin{proof}
Let $s: \cO_X \ra L$ be a non-zero section. It lifts to $E'$ if and only if the extension class of \eqref{eq3.1} is 
contained in the kernel of the map $H^1(L^* \otimes F) \stackrel{s^*}{\ra} H^1(F)$. We have to show that $\ker s^* \neq 0$.

The exact sequence $0 \ra L^* \stackrel{s^*}{\ra} \cO_X \ra \tau \ra 0$ tensorized with $F$ gives the long exact sequence
$$
0 \ra H^0(F) \ra H^0(F \otimes \tau) \ra H^1(L^* \otimes F) \stackrel{s^*}{\ra} H^1(F) \ra 0.
$$ 
Now $h^0(F \otimes \tau) = d_F + 1$ and according to Clifford's Theorem (see \cite[Theorem 2.1]{bgn}),
$$
h^0(F) \leq \frac{d_F + 3}{2},
$$
since $F$ is stable and $d_F$ is odd. Moreover $d_F > 1$. So $h^0(F) < h^0(F \otimes \tau)$ and 
hence $\ker s^*$ is non-zero. 
\end{proof}

\begin{theorem} \label{thm3.1}
Any semistable bundle of rank $3$ computing either $\gamma'_3$ or $\gamma_3$ is generically generated.
\end{theorem}

\begin{proof}
We give the proof for $E$ computing $\gamma'_3$. The proof for $\gamma_3$ is obtained by replacing 
$\gamma'_3$ by $\gamma_3$ throughout.
 
Let $E$ be semistable of rank 3 and degree $d$ computing $\gamma'_3$. 
Suppose there is a proper subbundle $F$ of $E$ with $h^0(F) = h^0(E) =3+s$. As usual we write $d_F = \deg F$.

If $F$ is not semistable, then $r_F = 2$ and we have an exact sequence \eqref{eq2.1} with $M$ and $N$ line bundles
with $d_M > d_N$. If $h^0(M) > \frac{s}{3} + 1$, then
\begin{equation}  \label{eq3.2}
\gamma'_3 \leq \gamma_1 \leq \gamma(M) < d_M -\frac{2s}{3} \leq \frac{d-2s}{3} = \gamma(E), 
\end{equation} 
a contradiction. Otherwise $h^0(N) \geq \frac{s}{3} + 1$ and $d_N < d_M \leq \frac{d}{3}$. Then we get the 
same contradiction replacing $M$ by $N$ in \eqref{eq3.2}. 

If $r_F = 1$, then, since $h^0(F) = 3+s$, we have $d_F \geq d_{s+2}$. Hence 
$$
d \geq 3d_{s+2} \geq 3(\gamma_1 + 2s + 4)
$$
by \cite[Lemma 4.6]{cl} (note that $d_{s+2} \leq \frac{d}{3} \leq g-1$). This gives 
$$
\gamma(E) \geq \gamma_1 + 2s + 4 - \frac{2s}{3} > \gamma_1 \geq \gamma'_3,
$$
a contradiction. 

So suppose $F$ is semistable of rank 2. If $d_F$ is even, consider 
$$
E' = F \oplus L
$$ 
with $L$ a line bundle of degree $d_L = \frac{d_F}{2}$ and $h^0(L) =1$. Then
$$
\gamma(E') < \gamma(E).
$$
This contradicts the assumption that $E$ computes $\gamma'_3$. 

If $d_F$ is odd, then by Lemmas \ref{lem3.1} and \ref{lem3.2} there exists a stable bundle $E'$ of rank 3 and degree
$\frac{3d_F + 1}{2}$ with $h^0(E') \geq h^0(E) + 1$. Hence
$\gamma(E') < \gamma(E)$ contradicting the 
minimality of $\gamma(E)$.  
\end{proof}

The following lemma shows that, if $\frac{d_3}{3} > \frac{d_2}{2}$, the conclusion of \cite[Corollary 4.12]{cl} may not hold.

\begin{lem} \label{rem3.4}
If $\frac{d_3}{3} > \frac{d_2}{2}$ and $d_2$ is even, then there exists a semistable bundle $E$ 
of rank $3$ with $h^0(E) = 4$ and $d_E = \frac{3d_2}{2} < d_3$.
\end{lem}

\begin{proof}
Note that $d_2 < 2d_1$, since otherwise $\frac{d_3}{3} > \frac{d_2}{2} = d_1$, a contradiction. 

Let $M$ be a line bundle of degree $d_2$ with $h^0(M) = 3$ and define $E_M$ by
$$
0 \ra E_M^* \ra H^0(M) \otimes \cO_C \ra M \ra 0.
$$
Then $E_M$ has rank 2, degree $d_2$ and $h^0(E_M) \geq 3$. In fact, $E_M$ is stable by 
\cite[Proposition 4.9(e)]{cl} and $h^0(E_M) = 3$ by \cite[Theorem 4.15(a)]{cl}. Now define 
$$
E = E_M \oplus N
$$
where $N$ is a line bundle of degree $\frac{d_2}{2}$ with $h^0(N) = 1$.
\end{proof}

\begin{prop} \label{prop3.5}
If $\frac{d_3}{3} > \frac{d_2}{2}$, then there exists a semistable bundle $E$ of rank $3$ with 
$h^0(E) = 4$ and 
\begin{equation} \label{eq3.3}
\frac{1}{2}(d_2-2) < \gamma(E) = \frac{1}{3} \left( \left[ \frac{3d_2 + 1}{2} \right] - 2 \right) \leq \frac{1}{3}(d_3 - 2).
\end{equation}
\end{prop}

\begin{proof}
If $d_2$ is even, we take $E$ to be one of the bundles described in Lemma \ref{rem3.4}. 
We have $\gamma(E) = \frac{1}{3}(\frac{3d_2}{2} - 2)$ as required.
The inequalities in \eqref{eq3.3} are straightforward computations.

If $d_2$ is odd, we take $F = E_M$ and a line bundle $L$ of degree $\frac{d_2 + 1}{2}$with $h^0(L) = 1$.
Lemma \ref{lem3.2} yields a bundle $E = E'$ with $h^0(E) =4$ 
and $d_E = \frac{3d_2+1}{2}$. By Lemma \ref{lem3.1}, $E$ is stable. \eqref{eq3.3} is again obvious.
\end{proof}

\begin{rem}
{\em Note that the bundles $E$ constructed in Proposition \ref{prop3.5} are not generated and that
$\gamma(E) > \frac{1}{2}(d_2 - 2) \geq \gamma_2$ by \cite[Corollary 5.3]{cl}. So, if $E$ computes 
$\gamma_3$, we must have $\gamma_3 > \gamma_2$; this is also an immediate consequence of Theorem \ref{thm2.4}.} 
\end{rem}

\begin{cor} \label{cor3.6}
Let $C$ be a smooth plane curve of degree $\delta \geq 7$. Then there exists a semistable bundle $E$ of rank $3$ 
and degree $d_E < d_3$ with $h^0(E) = 4$. If 
\begin{equation} \label{eq3.4}
\gamma'_3 \geq \frac{1}{3} \left( \left[ \frac{3 \delta + 1}{2} \right] - 2 \right),
\end{equation}
then $E$ computes $\gamma_3$. 
\end{cor}

\begin{proof}
For a smooth plane curve we have $d_2 = \delta$ and the hypotheses of Proposition \ref{prop3.5} apply.
If \eqref{eq3.4} holds, then it follows from \cite[Proposition 8.3]{cl} that 
$\gamma_3 = \frac{1}{3}([\frac{3 \delta + 1}{2}] -2) = \gamma(E)$, where $E$ is any of the bundles 
constructed in Proposition \ref{prop3.5}.
Moreover, 
$$
d_E = \left[\frac{3\delta + 1}{2}\right] < 2 \delta -2 = d_3.
$$ 
\end{proof}

\begin{rem}
{\em In a forthcoming paper \cite[Proposition 4.6 and Corollary 4.8]{lb}, we show that \eqref{eq3.4} holds for $\delta\ge10$. Thus, for such curves, there exists a non-generated bundle computing $\gamma_3$ and $\gamma_3>\gamma_2$.}
\end{rem}


\begin{thebibliography}{CAV}
\bibitem{b} E. Ballico:
\emph{Spanned vector bundles on algebraic curves and linear series}.
Rend. Istit. Mat. Univ. Trieste 27 (1995), no. 1-2, 137--156 (1996).
\bibitem{bgn} L. Brambila-Paz, I. Grzegorczyk, P. E. Newstead:
\emph{Geography of Brill-Noether loci for small slopes}.
J. Alg. Geom. 6 (1997), 645-669.

\bibitem{ckm} M. Coppens, C. Keem, G. Martens:
\emph{Primitive linear series on curves}.
Manuscr. Math. 77 (1992), 237 - 264.
\bibitem{cl} H. Lange and P. E. Newstead: 
\emph{Clifford indices for vector bundles on curves}.
To appear in: A. Schmitt (Ed.) Affine Flag Manifolds and Principal Bundles. Birkh\"auser.
\bibitem{lb} H. Lange and P. E. Newstead:
\emph{Lower bounds for Clifford indices in rank three}.
arXiv 0912.2618v1.

\end{thebibliography}
\end{document}